\documentclass[12pt]{article}
\usepackage{amsmath}
\usepackage{amssymb}
\usepackage{amsthm}
\newtheorem{theorem}{Theorem}
\newtheorem{lemma}[theorem]{Lemma}
\newtheorem{corollary}[theorem]{Corollary}
\newtheorem{proposition}[theorem]{Proposition}

\newtheorem*{thmM}{Main Theorem}
\font\fppl=pplr7t

\begin{document}

\def\kk{{\bf k}}
\def\kkk#1{{\bf k}\langle #1\rangle}
\def\phi{\varphi}
\def\leq{\leqslant}
\def\geq{\geqslant}

\fppl

\title{ An analogue of the Magnus problem for associative
algebras}

\author{\fppl V.~Dotsenko, N.~Iyudu and D.~Korytin}

\date{}

\maketitle

\begin{abstract} \fppl We prove an analogue of the Magnus theorem for
associative algebras without unity over arbitrary fields. Namely, if
an algebra is given by $n+k$ generators and $k$ relations and has an
$n$-element system of generators, then this algebra is a free
algebra of rank $n$.
\end{abstract}

Let we recall two group-theoretic statements concerning free
subgroups in finitely presented groups. The Magnus theorem
\cite{mag1,mag2} says that if a group $G$ is given by $n+k$
generators and $k$ relations and has a system of $n$ generators,
then $G$ is a free group of rank $n$. The following theorem was
proved by Romanowskii \cite{rom}: let $G$ be a group given by $n+k$
generators and $k$ relations. Then there exists a free subgroup of
$G$ of rank $n$. Moreover, generators of this free subgroup can be
chosen from the initial generators of the group. Analogues of both
theorems for semigroups were proved by Shneerson \cite{shne1}. The
second theorem (Lyndon condition) for associative algebras over
finite fields was announced in \cite{shne2} and proved in
\cite{khar}. We present here the proof of the analog of Magnus
theorem for associative algebras without unity over arbitrary field.

\begin{thmM}  Let $A$ be an associative algebra without unity given by
$n+k$ generators and $k$ relations over a field $\kk$. If $A$ has an
$n$-element system of generators $g_1,\dots,g_n$, then $A$ is a free
algebra of rank $n$.
\end{thmM}

Let us fix the presentation of $A$:
$$
A=\kkk{X}/I,\quad I={\rm id}(h_1,\dots,h_k),\quad h_i\in\kkk{X},
$$
where $X=\{x_1,\dots,x_{n+k}\}$. We assume that the free associative
algebra $\kkk{X}$ is supplied by the ordinary  degree function, with
all variables having degree 1. Denote the linear part (homogeneous
part of degree 1) of $g\in\kkk{X}$ by $Lg$. Let we mention two easy
general facts.

\begin{proposition} \label{p1}  Let $V=\{v_1,\dots,v_n\}$ be
a system of elements  $V\subset \kkk{X}$  and $W=\{w_1,\dots,w_n\}$
obtained from $V$ by a non-degenerate linear transformation. Then
algebraic independence of $V$ is equivalent to algebraic
independence of $W$.
\end{proposition}

\begin{corollary} \label{c1}  Let $V=\{v_1,\dots,v_n\}$, $V\subset \kkk{X}$ be
a system of elements of degree $1$. Then $V$ is algebraically
independent if and only if $V$ is linearly independent.
\end{corollary}

Let $g_1,\dots,g_n$ be a system of generators of $A$. For
convenience,  we denote $h_j$ by $g_{n+j}$. We also denote linear
parts of polynomials $g_1,\dots,g_{n+k}$ as $y_1,\dots,y_{n+k}$:
$Lg_j=y_j$.

\begin{theorem} \label{t1}   The system $\{y_1,\dots,y_{n+k}\}$ is
algebraically independent. \end{theorem}

\begin{proof}  \fppl According to Corollary~\ref{c1}, it suffices to verify
that $y_1,\dots,y_{n+k}$ are linearly independent. Since
$g_1,\dots,g_n$ generate the quotient
$A=\kkk{x_1,\dots,x_{n+k}}/{\rm id}(g_{n+1},\dots,g_{n+k})$, we have
that
\begin{equation}\label{e1}
x_i=\Phi_i(g_1,\dots,g_n)+d_i,\quad d_i\in{\rm
id}(g_{n+1},\dots,g_{n+k})
\end{equation}
for any $i=1,\dots,n+k$. Comparing the linear parts in these
relations, we obtain
\begin{equation}\label{e2}
x_i=\sum_{r=1}^{n+k} \alpha_r^iy_r,\quad \alpha_r^i\in\kk
\end{equation}
for any $i=1,\dots,n+k$. Hence the $(n+k)$-element system
$y_1,\dots,y_{n+k}$ is linearly independent.

\end{proof}

Consider the isomorphism $\phi$ of free algebras
$\phi:\kkk{x_1,\dots,x_{n+k}}\to \kkk{y_1,\dots,y_{n+k}}$ defined by
the formula $\phi(x_i)=\sum\limits_{r=1}^{n+k} \alpha_r^iy_r$.

\begin{lemma}\label{t2}  The system $g_1,\dots,g_n$ is
algebraically independent in the algebra
$$
\kkk{x_1,\dots,x_{n+k}}/{\rm id}(g_{n+1},\dots,g_{n+k})
$$
if and only if the system $\phi g_1,\dots,\phi g_n$ is algebraically
independent in the algebra
$$
\kkk{y_1,\dots,y_{n+k}}/{\rm id}(\phi g_{n+1},\dots,\phi g_{n+k}).
$$
\end{lemma}

\begin{proof}  \fppl Since $\phi$ is an isomorphism of free associative
algebras, the lemma follows.
\end{proof}

\begin{theorem}\label{t3} The system $\phi g_1,\dots,\phi g_n$ is algebraically
independent in the algebra
$$
A=\kkk{y_1,\dots,y_{n+k}}/{\rm id}(\phi g_{n+1},\dots,\phi g_{n+k}).
$$
\end{theorem}

\begin{proof}  \fppl By the described non-degenerate change of variables,
we ensure that:
\begin{equation} \label{e3}
L\phi g_i=y_i,\quad i=1,\dots,n+k.
\end{equation}
Indeed, since $\phi(x_i)$ is equal to $\sum\alpha_k^iy_k$ and $y_k$
are linearly independent, there exists an inverse transformation
$y_j=\sum\beta_k^jx_k$. Hence
$$
L\phi g_j=\phi Lg_j=\phi y_j=\sum_k\beta_k^j\phi
x_k=\sum_{k,i}\beta_k^j\alpha^k_i\phi y_i=y_j.
$$
Suppose there exists $\Phi\in\kkk{z_1,\dots,z_n}$, $\Phi\neq 0$ such
that
$$
\Phi(\phi g_1,\dots,\phi g_n)=d\quad \text{for}\quad d\in{\rm
id}(\phi g_{n+1},\dots,\phi g_{n+k}).
$$
The property (\ref{e3}) allows us to compare the minimal homogeneous
parts of $\Phi(\phi g_1,\dots,\phi g_n)$ and $d\in{\rm id}(\phi
g_{n+1},\dots,\phi g_{n+k})$. Homogeneous part of minimal degree of
the polynomial $\Phi(\phi g_1,\dots,\phi g_n)$ is
$M(y_1,\dots,y_n)=\sum\limits_{(i)}\gamma_{(i)}y_{i_1}\dots
y_{i_l}$, where $M(z_1,\dots,z_n)$ is the homogeneous part of
minimal degree of $\Phi(z_1,\dots,z_n)$. Indeed,
$M(y_1,\dots,y_n)\neq 0$ because $y_1,\dots,y_n$ are algebraically
independent. Hence $\Phi(\phi g_1,\dots,\Phi g_n)$ contains only
variables $y_1,\dots,y_n$ and does not contain
$y_{n+1},\dots,y_{n+k}$.

Now we shall prove a lemma, which shows that some minimal degree
monomials of all non-zero elements of the ideal ${\rm id}(\phi
g_{n+1},\dots,\phi g_{n+k})$ contain one of the variables
$y_{n+1},\dots,y_{n+k}$. This contradiction will show that
$\Phi(z_1,\dots,z_n)=0$ and therefore the elements $\phi g_1,\dots,
\phi g_n$ are algebraically independent.

Let us consider the free associative algebra
$\kkk{x_1,\dots,x_{n+k}}$ with the ordering on the variables
$x_1<{\dots}<x_{n+k}$. It induces the degree-lexicographical
ordering on the set of monomials $\langle X\rangle$: $x_{i_1}\dots
x_{i_k}<x_{j_1}\dots x_{j_m}$ if $k<m$ or $k=m$ and there exists $l$
such that $i_l<j_l$ and $i_s=j_s$ for $s<l$.

For $s\in \kkk{X}$ by $m(s)$ we denote the minimal monomial (with
respect to the above ordering) which appears in $s$ with non-zero
coefficient.

\begin{lemma}\label{l1} Let
$f_1,\dots,f_k\in\kkk{x_1,\dots,x_{n+k}}$ and the linear part of
$f_i$ is $x_i$ for any $i$: $Lf_i=x_i$. Let also $I={\rm
id}(f_1,\dots,f_k)$ be the ideal generated by $f_1,\dots,f_k$. Then
for any $s\in I$, its minimal monomial $m(s)$ contains at least one
of the variables $x_1,\dots,x_k$.
\end{lemma}

\begin{proof}  \fppl For an arbitrary presentation
$$
s=\Pi(s)=\sum_{i,\alpha_i} c_{\alpha_i}p_{\alpha_i}f_iq_{\alpha_i},
$$
where $c_{\alpha_i}\in \kk\setminus\{0\}$ and
$p_{\alpha_i},q_{\alpha_i}\in \langle X\rangle$, we define the {\it
parameter $\tau=\tau(\Pi(s))$ of the presentation $\Pi(s)$} as the
minimal monomial appearing in the presentation:
$$
\tau=\min_{i,\alpha_i,j}p_{\alpha_i}u_j(f_i)q_{\alpha_i},
$$
where $u_j(f_i)$ are monomials appearing in the polynomial $f_i$. We
shall prove that under our conditions for each $s\in {\rm
id}(f_1,\dots,f_k)$, there exists a presentation $\Pi(s)$ of $s$
with the parameter $m(s)$. We start with an arbitrary presentation
of $s$:
$$
s=\sum_{\alpha_i,\ i=1,\dots,k}
c_{\alpha_i}p_{\alpha_i}f_iq_{\alpha_i},
$$
where $c_{\alpha_i}\in \kk\setminus\{0\}$ and
$p_{\alpha_i},q_{\alpha_i}$ are monomials. Consider those elements
$p_{\alpha_i}f_iq_{\alpha_i}$ of this sum for which
$p_{\alpha_i}x_iq_{\alpha_i}$  equals to the parameter $\tau$ of
this presentation. Clearly $\tau\leq m(s)$. If $\tau=m(s)$, we are
through. Assume that $\tau<m(s)$ and consider
$$
M=\sum_{i,\alpha_i:p_{\alpha_i}x_iq_{\alpha_i}=\tau}
c_{\alpha_i}p_{\alpha_i}f_iq_{\alpha_i}.
$$
Since $\tau<m(s)$, we have $M=0$. We will obtain a new presentation
of $s$ if in $\Pi(s)$ we replace all $x_i$  with $1\leq i \geq k$ by
$f_i-(f_i-x_i)$ in all monomials $p_{\alpha_i}$, $q_{\alpha_i}$ from
the sum $M$. Since $M=0$, the same sum with $x_i$ replaced by $f_i$
also equals zero. Moreover, if we eliminate the last sum from the
new presentation, we again obtain a presentation of $s$ whose
parameter is greater than $\tau$ since $m(f_i-x_i)>x_i$. Hence we
have obtained a presentation of $s$ with a greater parameter. If the
new parameter is still less that $m(s)$, we can apply the same
procedure to get another presentation with greater parameter etc.
Since it can not become greater than $m(s)$, the process will
terminate and we will obtain a presentation of $s$ with parameter
$m(s)$. Since the parameter of any presentation via $f_1,\dots,f_k$
contains at least one of the variables $x_1,\dots,x_k$, we obtain
the statement of the lemma.
\end{proof}

As it was already mentioned above, the proof of Lemma~\ref{l1}
completes the proof of Theorem~\ref{t3}.
\end{proof}

Now the Main Theorem follows from Lemma~\ref{t2} and
Theorem~\ref{t3}.

\smallskip

{\bf Acknowledgements.} \ We would like to thank Professor Viktor
Latyshev, who draw our attention to this very nice combinatorial
problem.

\small\rm

\end{document}